\documentclass[12pt,a4paper]{article}

\usepackage[cp1251]{inputenc}
\usepackage[english,russian]{babel}

\usepackage{amsmath}
\usepackage{amssymb}
\usepackage{amsfonts}
\usepackage{amsthm}
\usepackage[matrix,arrow]{xy}
\usepackage{longtable}
\usepackage{graphicx}
\usepackage{amscd}
\usepackage{pb-diagram}

\newtheorem{theorem}{Theorem}[section]
\newtheorem{stm}[theorem]{Утверждение}

\newtheorem{lemma}[theorem]{Lemma}
\newtheorem{imp}[theorem]{Corollary}

\theoremstyle{definition}

\newcommand{\mb}{\mathbb}

\newcommand{\Mo}{\overline{\mathcal{M}}_{g;m}}
\newcommand{\Hu}{P\overline{\mathcal{H}}_{g;\kappa}}
\newcommand{\Huo}{P\overline{\mathcal{H}}_{0;\kappa}}

\renewcommand{\proof}{{\it Proof.}\;\;}

\title{Degrees of cohomological classes of multisinguliarities in Hurwitz spaces of rational functions}
\author{B.\;Bychkov \footnote{National Research University "Higher School of Economics", Russian Federation. This work has been funded by the  Russian Academic Excellence Project "5-100". The research is also supported by the Russian Science Foundation grant, project 16-11-10316.}}
\date{}

\begin{document}
\maketitle

Hurwitz numbers enumerate ramified coverings of the $2$-dimensional sphere with prescribed ramifications. They play one of the central roles in algebraic geometry and theory of integrable systems during the last two decades. All known methods for computing Hurwitz numbers allow us to receive explicit formulae for families of these numbers only in few specific cases.
One of these methods is based on the theory of universal polynomials. It was developed by M.~Kazarian and S.~Lando, and allow one to compute degrees of certain cohomology classes in the Hurwitz spaces.

Hurwitz spaces are spaces of meromorphic functions with prescribed orders of poles on curves of a given genus. In this paper we receive new formulae for degrees of strata of genus~$0$ Hurwitz spaces. These strata correspond to functions with two non simple critical values with prescribed partitions of multiplicities of preimages. Herewith one of the preimages has an arbitrary multiplicity and another has a multiplicity of codimension 1.
Universal cohomological expressions for strata of codimension~$1$ in Hurwitz spaces were received by M.~Kazarian and S.~Lando in \cite{KL}. We get new explicit formulae for certain families of genus zero double Hurwitz numbers. These formulae were conjectured by M.~Kazarian on the base of computer experiments. We have all prerequisites to obtain more general explicit formulae for Hurwitz numbers and we think that our methods will let to prove more general identities on series of functions of partitions, which encode relations in cohomology rings of Hurwitz spaces.

The paper splits into two parts. In the first part we give a description of Hurwitz spaces, introduce all notions and represent main results. In the second one we collect all proofs, and at the end of the paper we give a proof of a combinatorial {\it Kazarian identity}.

We are grateful to S.~Lando who was the initiator of research in this direction; as well, he is the author of the proof of lemma \ref{l42}; and to M.~Kazarian who formulated the conjecture, on which we underpin the paper, and gave a detailed explanation of the paper \cite{KL}; and to D.~Ilinskiy for useful discussions.

\section{Hurwitz spaces}
\subsection{Stratification of the discriminant of Hurwitz spaces}

Let $\mathcal{H}_{g;\kappa}$, $\kappa=(k_1,\dots,k_n)$, be the space of meromorphic functions of degree $k_1+\ldots+k_m=n$ on genus $g$ algebraic curves with the following properties:
\begin{itemize}
\item[$-$] every function has $m$ numbered poles of given orders $k_1,\ldots,k_m$;
\item[$-$] the sum of all critical values of the function is equal to zero.
\end{itemize}

From the topological point of view, every such function is an $n$-sheeted ramified covering of the
$2$-dimensional sphere by a genus $g$ surface.

According to \cite{ELSV}, this space is a smooth complex orbifold (moreover, for $g=0$ or sufficiently large $n$, $\mathcal{H}_{g;\kappa}$ is a complex manifold). Let $\mathcal{M}_{g;m}$ be the moduli space of complex curves of genus $g$ with $m$ marked points, then $\mathcal{H}_{g;\kappa}$ is fibered over $\mathcal{M}_{g;m}$: one can associate to the function its domain with the $m$ marked poles.

Let us denote by $\overline{\mathcal{H}}_{g;\kappa}$ the completion of the space $\mathcal{H}_{g;\kappa}$ that consists of {\it stable} meromorphic functions (see \cite{ELSV}, \cite{La}, \cite{Fu}).
Stable maps are defined on nodal curves (curves whose only singuliarities are simple self intersections) with $m$ marked smooth points. A map is stable if it has finitely many automorphisms. The boundary $\overline{\mathcal{H}}_{g;\kappa} \setminus \mathcal{H}_{g;\kappa}$ of the completed Hurwitz space consists of stable functions on singular curves.

Let us denote by $\Mo$ the moduli space of {\it stable} curves which is the Deligne--Mamford compactification of the moduli space of smooth curves $\mathcal{M}_{g;m}$. Recall that a nodal curve with smooth marked points is said to be stable if its automorphism group is of finite order.
The projection $\mathcal{H}_{g;\kappa} \rightarrow \mathcal{M}_{g;m}$ can be extended to the projection $\overline{\mathcal{H}}_{g;\kappa} \rightarrow \Mo$. The fiberwise projectivization  $\Hu$ of the bundle $\overline{\mathcal{H}}_{g;\kappa}$ is a compact complex orbifold.

The main subject of our investigation are submanifolds in Hurwitz spaces consisting of functions with degenerate finite critical values. A critical value in the image is said to be {\it nondegenerate} if it has $n-1$ distinct preimages: one of the preimages has multiplicity~$2$, while other $n-2$ preimages are smooth points of the covering, and it is said to be {\it degenerate} otherwise.
By the Riemann--Hurwitz formula, a generic meromorphic function of degree $n$ on a genus $g$ curve with fixed multiplicities $k_1,\ldots,k_m$ of preimages of infinity has $n+m+2g-2$ nondegenerate ramification points. Functions with less than $n+m+2g-2$ critical values in the image form the {\it discriminant} in $\Hu$. To each critical value in the image we associate a partition $\mu$ of $n$, which is the unordered tuple of the multiplicities of the preimages of this point. Let us denote by $P\overline{\mathcal{H}}_{g;\kappa;\mu_1;\ldots;\mu_l}$ the closure in $\Hu$ of the set of functions having prescribed ramifications; the index in the notation consists of the genus $g$, the multiplicities of the poles and the tuple of the partitions of the multiplicities of preimages over finite degenerate critical values. This is a {\it discriminant stratum}.

Each discriminant stratum $P\overline{\mathcal{H}}_{g;\kappa;\mu_1;\ldots;\mu_l}$ is a complex submanifold of pure dimension in $\Hu$ and therefore it determines, by Poincar\'e duality, a homogeneous element in the cohomology ring $H^*(\Hu)$. Let us denote this element by $\sigma_{g;\kappa;\mu_1;\ldots;\mu_l},$
\begin{equation*}
\sigma_{g;\kappa;\mu_1;\ldots;\mu_l} = [P\overline{\mathcal{H}}_{g;\kappa;\mu_1;\ldots;\mu_l}] \in H^*(\Hu).
\end{equation*}

The {\it degree} of a class $\sigma_{g;\kappa;\mu_1;\ldots;\mu_l}$ is the intersection index of $\sigma_{g;\kappa;\mu_1;\ldots;\mu_l}$ with the complementary degree of the first Chern class of the tautological line bundle over $P\overline{\mathcal{H}}_{g;\kappa;\mu_1;\ldots;\mu_l}$. More generally, let $A$ be a complex manifold and let the multiplicative group $\mb{C}^*$ of nonzero complex numbers act on~$A$ without fixed points. Let $B$ be the manifold of orbits of this action and suppose $B$ is compact. Let $\Psi = c_1(\mathcal{O}(1)) \in H^2(B)$ be the tautological class of this action. The {\it degree} $\mathrm{deg}\;\beta$ of an arbitrary class $\beta \in H^*(B)$ is the result of the pairing:
\begin{equation*}
\mathrm{deg}\;\beta = \int\limits_{B} \frac{\beta}{1 - \Psi} = \int\limits_{B} (1+\Psi+\Psi^2+\ldots)\beta.
\end{equation*}

\subsection{Main results}

In what follows, we assume that $g=0$. Discriminant strata of codimension 1 in Hurwitz spaces correspond to two types of degeneration of functions: the function may have a critical value with a triple preimage (the {\it caustic}) or may have a critical value with two double preimages (the {\it Maxwell stratum}).
These are, respectively, the strata $P\overline{\mathcal{H}}_{0;\kappa;1^{n-3}3^1}$ and $P\overline{\mathcal{H}}_{0;\kappa;1^{n-4}2^2}$.

\begin{theorem}
The following equalities hold:
\begin{multline*}
\mathrm{deg}\, \sigma_{0;\kappa;1^{n-3}3^1} =\\
 =n^{m-4} \prod\limits_{i=1}^m \frac{k_i^{k_i}}{k_i!}\left(n^2+n\left(3(m-2) - \left(\frac{1}{k_1} + \ldots + \frac{1}{k_m}\right)\right) + \frac{3}{2}(m-2)(m-3) \right).
\end{multline*}
\begin{multline*}
\mathrm{deg}\, \sigma_{0;\kappa;1^{n-4}2^2} =\\
 =n^{m-4} \prod\limits_{i=1}^m \frac{k_i^{k_i}}{k_i!}\left(-2n^2+n\left(m^2 - 5(m-2) + \left(\frac{1}{k_1} + \ldots + \frac{1}{k_m}\right)\right) - 2(m-2)(m-3) \right).
\end{multline*}
\label{T12}
\end{theorem}

The statement of Theorem~\ref{T12} is a generalization of results from~\cite{KL}, where the degrees of the discriminant strata with two degenerate critical values were calculated for all strata of codimension $2$:
\begin{eqnarray*}
\mathrm{deg}\,\sigma_{1^{n-3}3^1;1^{n-3}3^1} &=&\frac{3}{8}n^{n-6}(27n^2 - 137n+180),
\\
\mathrm{deg}\,\sigma_{1^{n-3}3^1;1^{n-4}2^2} &=& 3n^{n-6}(n-3)(3n^2-15n+20),
\\
\mathrm{deg}\,\sigma_{1^{n-4}2^2;1^{n-4}2^2} &=& 4n^{n-6}(2n^3-16n^2+43n-40).
\end{eqnarray*}

For example, if we substitute $\kappa=1^{n-3}3^1$, i.e., $k_1=\dots=k_{n-3}=1$, $k_{n-2}=3$ and $m=n-2$ into the first statement of the theorem, then we obtain the first of these equalities. The factor $\frac{1}{2}$ appears because we don't fix the order of the critical values.

The plan of the proof of Theorem~\ref{T12} is like follows: firstly, we reduce the computation of the degree of the stratum on the left hand side to computation of the degrees of some basic cohomological classes in $\Huo$,
and then we compute the latter degrees.

The first step of the proof is based on paper~\cite{KL}, where Kazarian and Lando deduced universal formulae for the
cohomology classes of codimension one strata in terms of certain basic classes in arbitrary generic families of rational functions:
\begin{equation}
\sigma_{0;1^{n-3}3^1} = -p_*\zeta + 3p_*\psi + 2p_*\psi^2 - \delta_{0,0}.
\label{sigma1}
\end{equation}

\begin{equation}
\sigma_{0;1^{n-4}2^2} = \frac{1}{2}\xi_0^2 -2p_*\zeta - 5p_*\psi - 3p_*\psi^2 + \delta_{0,0}.
\label{sigma2}
\end{equation}

By virtue of universality, these formulae hold in the Hurwitz spaces $\Hu$ as well. We will define the cohomology classes
entering the right hand side of (\ref{sigma1}) and (\ref{sigma2}) in Section \S~2.

\subsection{Hurwitz numbers}

The geometry of discriminants in Hurwitz spaces is closely related to the problem of computing Hurwitz numbers.

Two ramified coverings $f_1:C_1\to S^2$ and $f_2:C_2\to S^2$ of the 2-dimensional sphere are said to be {\it isomorphic\/} if there is a homeomorphism $h:C_1\to C_2$ such that $f_1=f_2\cdot h.$ As we discussed above every critical value in the image compares with a partition $\mu$ of multiplicities of it preimages. Let us define the {\it Hurwitz number} $h_{g;\mu_1;\mu_2;\ldots}$ by the following equality:
\begin{equation*}
h_{g;\mu_1;\mu_2;\ldots} = \sum\limits_f \frac{1}{|\mathrm{Aut}(f)|},
\end{equation*}
where the summation is over all isomorphism classes of $n$-sheeted ramified coverings $f$ of the 2-dimensional sphere by a genus $g$ surface with partitions $\mu_1,\mu_2,\ldots$ over degenerate critical values and by $|\mathrm{Aut}(f)|$ we denote the order of the automorphism group of this covering. The {\it Hurwitz problem} is the problem about computation of Hurwitz numbers. The following theorem  \cite{ELSV,KL} relates the Hurwitz problem to that of computing the degrees of discriminant strata in Hurwitz spaces.

\begin{theorem}
The following equation holds
\begin{equation*}
h_{g;\mu_1;\mu_2;\ldots} = \frac{|\mathrm{Aut(\mu_1,\mu_2,\ldots)|}r!}{n!} \mathrm{deg}\,P\overline{\mathcal{H}}_{g;\mu_1;\mu_2,\ldots},
\end{equation*}
here $\mathrm{Aut(\mu_1,\mu_2,\ldots)}$ denotes the product of factorials of the numbers of coinciding partitions and~$r$ denotes the total number of nondegenerate critical values of functions from $P\overline{\mathcal{H}}_{g;\mu_1;\mu_2,\ldots}$.
\label{T22}
\end{theorem}

{\it Double Hurwitz numbers\/} are the numbers $h_{g;\mu_1;\mu_2}$. Unlike simple Hurwitz numbers $h_{g;\mu}$ (see \cite{ELSV}) we don't know any useful universal closed formula for double Hurwitz numbers.

Let us briefly describe the current state of the study of double Hurwitz numbers.

In paper \cite{GJV}, explicit formulae for 2-part and 3-part genus $0$ Hurwitz numbers were obtained: 
this is the case where the preimage of one of the critical values consists of either two or three points, respectively. 
However, these formulae are hard to apply, since they contain a summation over all possible decompositions
 of a partition into a sum of two (or three) partitions. In the same paper, Goulden, Jackson and Vakil showed that the numbers $h_{g;\mu_1;\mu_2}$ are polynomials in the  parts of the partitions $\mu_1$ and $\mu_2$.

In paper \cite{SSV}, an explicit form of these polynomial for genus $0$ double Hurwitz numbers
inside so called chambers was deduced. In paper \cite{CJM}, polynomiality of double Hurwitz 
numbers inside chambers for arbitrary genus was proved. Let us note that in our case every such chamber 
corresponds to exactly one Hurwitz number.

A.~Okounkov in paper \cite{Ok} showed that the exponent of a generating series for double Hurwitz numbers
is a $\tau$-function of the Toda hierarchy. This property produces recurrent relations for computing double 
Hurwitz numbers of arbitrary genus, but it does not lead to explicit formulae. Topological recursion methods \cite{MSS} 
also provide certain recurrent relations for double Hurwitz numbers.

Theorems \ref{T12} and~\ref{T22} imply the following corollaries:

\begin{imp}
The following equalities hold
\begin{multline*}
h_{0;\kappa;1^{n-3}3^1}  =\frac{|\mathrm{Aut(\kappa,1^{n-3}3^1)|(n+m-4)!}}{n!}\times \\
\times n^{m-4} \prod\limits_{i=1}^m \frac{k_i^{k_i}}{k_i!}\left(n^2+n\left(3(m-2) - \left(\frac{1}{k_1} + \ldots + \frac{1}{k_m}\right)\right) + \frac{3}{2}(m-2)(m-3) \right).
\end{multline*}
\begin{multline*}
h_{0;\kappa;1^{n-4}2^2}  =\frac{|\mathrm{Aut(\kappa,1^{n-4}2^2)|(n+m-4)!}}{n!}\times \\
\times n^{m-4} \prod\limits_{i=1}^m \frac{k_i^{k_i}}{k_i!}\left(-2n^2+n\left(m^2 - 5(m-2) + \left(\frac{1}{k_1} + \ldots + \frac{1}{k_m}\right)\right) - 2(m-2)(m-3) \right).
\end{multline*}
\end{imp}

\section{Computations in cohomology rings of Hurwitz spaces}

Let us denote by $\overline{\mathcal{H}}_{0;\kappa,1}$ the space of stable rational functions of degree~$k_1+\ldots+k_m$ with orders of poles equal to $(k_1,\ldots,k_m)=\kappa$ at the first $m$ of the $m+1$ marked points and we require that the $(m+1)$st marked point is
a critical one. This space is a vector bundle over the moduli space $\overline{\mathcal{M}}_{0;m+1}$. 
It has the fiberwise projectivisation $P\overline{\mathcal{H}}_{0;\kappa,1}$. Let us consider the following commutative diagram

\begin{equation}\label{MCD}
\begin{CD}
P\overline{\mathcal{H}}_{0;\kappa,1} @ >\pi>> \overline{\mathcal{M}}_{0;m+1}\\
@V{p}VV @VV{}V\\
\Huo @>>> \overline{\mathcal{M}}_{0;m}
\end{CD}
\end{equation}
\noindent
The vertical arrows here denote the forgetful maps of the $(m+1)$st marked point 
and the horizontal ones are projection maps of the projectivized bundle to the base.

The map $p$ is proper since the $(m+1)$st marked point cannot coincide with any of the first $m$ marked points.

Since the space $P\overline{\mathcal{H}}_{0;\kappa,1}$ is a projectivization,
 we can consider a natural second cohomology class, namely, the first Chern class of the tautological line bundle
 over this space. We will denote it by $\zeta = c_1(\mathcal{O}(1))\in H^2(P\overline{\mathcal{H}}_{0;\kappa,1}).$

We also would like to consider the class  $\psi\in H^2(P\overline{\mathcal{H}}_{0;\kappa,1})$ equal to the first Chern class of the relative dualizing sheaf of the map $p: P\overline{\mathcal{H}}_{0;\kappa,1} \rightarrow \Huo$. Let us denote by $\Delta \in H^4(P\overline{\mathcal{H}}_{0;\kappa,1})$ the class represented by the set of singular points of the fibers of the map $p$. Let us introduce the classes $\delta_{k,l} = p_*(N^k\Delta^{l+1})\in H^{2k+4l+2}(\Huo)$, $k=0,1,2,\dots$, $l=0,1,2,\dots$. Here $N$ 
is the first Chern class of the normal bundle over the subvariety of singular points of fibers of the map $p$ see \cite{KL}.

Let us denote by $\mathrm{deg}_{\kappa}\;\beta$ the degree of the class $\beta\in H^*(\Huo)$ in the space $\Huo$. We will use the next two lemmas during the proof of the main theorem.

\begin{lemma} The following equalities hold
\begin{eqnarray*}
\mathrm{deg}_{k_1,\ldots,k_m}\;p_*\zeta^k&=& n^{m-2}\prod\limits_{i=1}^m \frac{k_i^{k_i}}{k_i!},\\
\mathrm{deg}_{k_1,\ldots,k_m}\;p_*\psi^k&=& \binom{m-2}{k}n^{m-2-k}\prod\limits_{i=1}^m \frac{k_i^{k_i}}{k_i!}.
\end{eqnarray*}
\label{L21}
\end{lemma}

\begin{lemma} The following equality holds
\begin{equation*}
\deg_{k_1,\ldots,k_m}\delta_{0,0} = \left(n\cdot\left(\frac{1}{k_1} + \ldots + \frac{1}{k_m}\right) - \frac{1}{2}(m-2)(m-3)\right)n^{m-4} \prod\limits_{i=1}^m \frac{k_i^{k_i}}{k_i!}.
\end{equation*}
\label{L23}
\end{lemma}

\subsection{Proof of Lemma \ref{L21}}

Let $\beta \in H^*(P\overline{\mathcal{H}}_{0;\kappa,1}) $, then
\begin{equation*}
\mathrm{deg}\;p_*(\beta) = \int\limits_{\Huo} \frac{p_*(\beta)}{1-p_*(\zeta)} = \int\limits_{P\overline{\mathcal{H}}_{0;\kappa,1}} \frac{\beta}{1-\zeta} =
\end{equation*}
\begin{equation}
 = \prod\limits_{i=1}^m \frac{k_i^{k_i}}{k_i!} \int\limits_{\overline{\mathcal{M}}_{0;m+1}} \frac{\pi_*(\beta)}{(1-k_1\psi_1)\ldots(1-k_m\psi_m)}.\label{pr1}
\end{equation}
The last equality is the result of computation of the total Segre class of the bundle $\Huo\to \overline{\mathcal{M}}_{0;m}$, see \cite{ELSV}.

Here by  $\psi_i$, $i=1,\dots,m+1$ we traditionally denote the first Chern classes of the line bundles $\mathcal{L}_i$ over $\overline{\mathcal{M}}_{0;m+1}$.
The fiber of the bundle $\mathcal{L}_i$ at a point $(C;x_1,\ldots,x_{m+1})$ of the moduli space is the cotangent line to the curve $C$ at the point $x_i$.
In our computations, we only will meet classes $\pi_*(\beta)$ that are equal to the class $\psi_{m+1}$ (i.e. the class $\psi$ at the $(m+1)$st marked point) in some degree.
Therefore, the computation of degrees of cohomological classes of strata of Hurwitz spaces reduces 
to computation of integrals of monomials of the $\psi$-classes over moduli spaces of curves. The results are known from the work of Witten~\cite{Wi}:

\begin{equation}
\int\limits_{\overline{\mathcal{M}}_{0;m}} \psi_1^{l_1}\ldots\psi_m^{l_m} = \binom{m-3}{l_1,\ldots,l_m}
\qquad\text{here }l_1+\dots+l_m=m-3. \label{pr2}
\end{equation}

We complete the proof of Lemma \ref{L21} by applying Eqs.~(\ref{pr1}) and (\ref{pr2}) to the classes $\zeta$ и $\psi$:

\begin{eqnarray*}
\mathrm{deg}_{\kappa}\;p_*\zeta^k &=& \prod\limits_{i=1}^m \frac{k_i^{k_i}}{k_i!} \int\limits_{\overline{\mathcal{M}}_{0;m+1}} \frac{1}{(1-k_1\psi_1)\ldots(1-k_m\psi_m)}\\
&=&\left(\sum_{l_1+\dots+l_m=m-2}k_1^{l_1}\dots k_m^{l_m}\int\limits_{\overline{\mathcal{M}}_{0;m+1}}\psi_1^{l_1}\dots
 \psi_m^{l_m}\right)\prod\limits_{i=1}^m \frac{k_i^{k_i}}{k_i!}\\
&=&\left(\sum_{l_1+\dots+l_m=m-2}\binom{m-2}{l_1,\dots,l_m}k_1^{l_1}\dots k_m^{l_m}\right)\prod\limits_{i=1}^m \frac{k_i^{k_i}}{k_i!}\\
&=& (k_1+\dots+k_m)^{m-2}\prod\limits_{i=1}^m \frac{k_i^{k_i}}{k_i!}= n^{m-2}\prod\limits_{i=1}^m \frac{k_i^{k_i}}{k_i!}.
\end{eqnarray*}
Similarly,
\begin{equation*}
\mathrm{deg}_{\kappa}\;p_*\psi^k = \prod\limits_{i=1}^m \frac{k_i^{k_i}}{k_i!} \int\limits_{\overline{\mathcal{M}}_{0;m+1}} \frac{\psi^k_{m+1}}{(1-k_1\psi_1)\ldots(1-k_m\psi_m)} =  \binom{m-2}{k}n^{m-2-k}\prod\limits_{i=1}^m \frac{k_i^{k_i}}{k_i!}.
\end{equation*}
$\blacktriangleleft$

\subsection{Proof of Lemma \ref{L23}}

To each partition, its Young diagram is associated in a standard way. 
Let $\kappa=(k_1,\ldots,k_m)$, $k_1+\ldots+k_m=n$, be a partition of~$n$,
$|\mathrm{Aut}(\kappa)|$ is the product of factorials of the numbers of coincides parts in the partition $\kappa$, $m=l(\kappa)$ 
is the length of a partition. Let us call a diagram $\mu$ a {\it subdiagram} of the diagram $\kappa$, $\mu<\kappa$, 
if $\mu$ is obtained from $\kappa$ by deleting some rows (parts). The {\it sum} of two diagram $\mu$ and $\lambda$ is the diagram $\mu\oplus\lambda$  obtained from the summands as the union of the sets of their parts.

\begin{stm}
The degree of the class $\delta_{0,0}$ is equal to
\begin{equation}
\deg_{k_1,\ldots,k_m}\delta_{0,0} = \frac{1}{2}  \sum\limits_{\mu\oplus\lambda=\kappa, \mu<\kappa} \frac{|\lambda|^{l(\lambda)-2}}{|\mathrm{Aut}(\lambda)|} \frac{|\mu|^{l(\mu)-2}}{|\mathrm{Aut}(\mu)|}|\mathrm{Aut}(\kappa)| \prod\limits_{i=1}^m \frac{k_i^{k_i}}{k_i!}
\label{St24}
\end{equation}
\label{delta}
\end{stm}
\proof
The class $\delta_{0,0}$ is the class dual to the subvariety of the singular points of fibers of the map $p$. 
A generic singular fiber is a singular rational curve with one singular point, i.e. a union of two smooth rational curves intersecting transversely at one point. If a smooth fiber degenerates into a generic singular fiber,
then the poles of the corresponding rational function are distributed between two irreducible components of the singular fiber.
Tuples of orders of poles form a partition $\lambda$ on one of the irreducible components and a partition $\mu$ on another irreducible component, where $\lambda\oplus\mu=\kappa$, $|\lambda|+|\mu|=n$.
The degree of the stratum $\delta_{0,0}$ is equal to the sum of the degrees of its components corresponding to 
various decompositions $\lambda,\mu$, whence we get Eq.~(\ref{St24}).
$\blacktriangleleft$

Note that each irreducible component of a generic singular fiber carries at least one pole.

To finish the proof of Lemma \ref{L23} we have to prove the following identity:

\begin{multline}
\frac{1}{2}  \sum\limits_{\mu\oplus\lambda=\kappa, \mu<\kappa}
 \frac{|\lambda|^{l(\lambda)-2}}{|\mathrm{Aut}(\lambda)|} \frac{|\mu|^{l(\mu)-2}}{|\mathrm{Aut}(\mu)|}|\mathrm{Aut}(\kappa)|= \\
=n^{m-4} \left(n\cdot\left(\frac{1}{k_1} + \ldots + \frac{1}{k_m}\right) - \frac{1}{2}(m-2)(m-3)\right).
\label{St25}
\end{multline}

M.~Kazarian conjectured this identity on the base of computer experiments. Below we would like to call this identity the {\it Kazarian combinatorial identity}. We will prove it in Secs. 2.4 and 2.5.

\subsection{Proof of the main theorem \ref{T12}}

Here we reproduce universal formulae for the cohomological classes Poincar\'e dual 
to codimension~$1$ discriminant strata from~\cite{KL}:

\begin{equation}
\sigma_{1^{n-3}3^1} = -p_*\zeta + 3p_*\psi + 2p_*\psi^2 - \delta_{0,0},
\label{proof1}
\end{equation}

\begin{equation}
\sigma_{1^{n-4}2^2} = \frac{1}{2}\xi_0^2 -2p_*\zeta - 5p_*\psi - 3p_*\psi^2 + \delta_{0,0}.
\label{proof2}
\end{equation}

Here by $\xi_0\in H^*(P\overline{\mathcal{H}}_{0;\kappa})$ we denote the class dual to the subvariety of critical points of functions in $P\overline{\mathcal{H}}_{0;\kappa}.$ The classes $\zeta,\; \psi$ and $\delta_{0,0}$ were introduced in \S~2.

To complete the proof, let us substitute results of Lemmas \ref{L21} and \ref{L23} into Eqs.~(\ref{proof1}) and (\ref{proof2}):

\begin{multline*}
\mathrm{deg}\, \sigma_{0;1^{n-3}3^1;\kappa} =\mathrm{deg}_{\kappa}\;\sigma_{1^{n-3}3^1} =\\
 =\prod\limits_{i=1}^m \frac{k_i^{k_i}}{k_i!} ( -n^{m-2} + 3\binom{m-2}{1}n^{m-3} + 2\binom{m-2}{2}n^{m-4} -\\
 - \left( n\left(\frac{1}{k_1} + \ldots + \frac{1}{k_m}\right)  -\binom{m-2}{2}\right)n^{m-4}  )=\\
= n^{m-4} \prod\limits_{i=1}^m \frac{k_i^{k_i}}{k_i!}\left(n^2+n\left(3(m-2) - \left(\frac{1}{k_1} + \ldots + \frac{1}{k_m}\right)\right) + \frac{3}{2}(m-2)(m-3) \right).
\end{multline*}

\begin{multline*}
\mathrm{deg}\, \sigma_{0;1^{n-4}2^2;\kappa} = \mathrm{deg}_{\kappa}\;\sigma_{1^{n-4}2^2} =\\
 =\prod\limits_{i=1}^m \frac{k_i^{k_i}}{k_i!} ( -n^{m-3}(2n-4)^2 -2n^{m-2}-5\binom{m-2}{1}n^{m-3} - 3\binom{m-2}{2}n^{m-4} +\\
 + \left(n\left(\frac{1}{k_1} + \ldots + \frac{1}{k_m}\right) - \binom{m-2}{2}\right)n^{m-4}  )=\\
 =n^{m-4} \prod\limits_{i=1}^m \frac{k_i^{k_i}}{k_i!}\left(-2n^2+n\left(m^2 - 5(m-2) + \left(\frac{1}{k_1} + \ldots + \frac{1}{k_m}\right)\right) - 2(m-2)(m-3) \right).
\end{multline*}

The theorem is proved.
$\blacktriangleleft$

\subsection{Kazarian identity}

In the next two sections we present a proof of the combinatorial identity~(\ref{St25}).

Let $\tau=(t_1,\ldots,t_m),\; t_1+\ldots+t_m=n$ be a partition.
Let us denote by $|\mathrm{Aut}(\tau)|$ the product of the factorials of the numbers of coinciding
parts in the partition $\tau$ and denote by $l(\tau)$ the length of~$\tau$.

Let $M$ be a finite set, $m=|M|$.
Denote by $t_M$ the sum of variables $t$ indexed by the set $M$,
$$
t_M=\sum_{i\in M} t_i.
$$

For each subset $I\subset M,\; I=\{j_1,\ldots,j_i\},$ denote by $\tau(I)$ the tuple $\{t_{j_1},\ldots,t_{j_i}\}$.

\begin{lemma}
Take two diagrams $\mu$ and $\lambda$ such that $\mu \oplus \lambda = \tau$;
 then the number of partitions of the set $M$ into two subsets $M = I \sqcup J$ such that $\tau(I) = \mu$, $\tau(J) = \lambda$ equals $\frac{|\mathrm{Aut}(\tau)|}{ |\mathrm{Aut}(\mu)| \cdot |\mathrm{Aut}(\lambda)|}$.
\label{L31}
\end{lemma}
\proof
The number of partitions of the set $M$ into two subsets with the above properties is grater than~$1$ if there are coinciding parts
among the parts $t_i$ of the partition $\tau$. The number of ways to pick some of these coinciding parts is equal to a binomial coefficient. The product of these binomial coefficients over all coinciding parts in the partition $\tau$ is equal to $\frac{|\mathrm{Aut}(\tau)|}{ |\mathrm{Aut}(\mu)| \cdot |\mathrm{Aut}(\lambda)|}$.
$\blacktriangleleft$

Now we are ready to transform the left hand side of the equality (\ref{St25}). Let us recall that $ |I|=l(\lambda),\; |J|=l(\mu),\; p=|\lambda| = t_I,\; q = |\mu|=t_J$. As a consequence of Lemma \ref{L31}, we get the following equation:

\begin{equation*}
\frac{1}{2}  \sum\limits_{\mu\oplus\lambda=\kappa, \mu<\kappa}\sum\limits_{\lambda\vdash p,\mu\vdash q} \frac{p^{l(\lambda)-2}}{|\mathrm{Aut}(\lambda)|} \frac{q^{l(\mu)-2}}{|\mathrm{Aut}(\mu)|}|\mathrm{Aut}(\kappa)| =
\frac{1}{2}\sum_{I \sqcup J= M} t_I^{|I|-2} t_J^{|J|-2}.
\end{equation*}

Rewrite the Kazarian identity in the form

\begin{equation}
\sum_{I \sqcup J = M} t_I^{|I|-2} t_J^{|J|-2} = t_M^{|M|-4} \left(2t_M\cdot\left(\frac{1}{t_1} + \ldots + \frac{1}{t_{|M|}}\right) - (|M|-2)(|M|-3)\right),\label{tk1}
\end{equation}
where summation is carried over all ordered partitions of the set~$M$ into two nonintersecting nonempty subsets.

\subsection{Abel polynomials}

We split the proof of identity (\ref{tk1}) into several steps.

First, notice that the left hand side and the right hand side of identity (\ref{tk1}) considered as functions
 in the variables $t_i$ are rational functions and have poles of order~$1$ at the points $t_i=0,$ $i\in M$. 
 The principal parts of the two functions at these points coincide. Indeed in the left hand side 
 a pole at the point $t_i=0$ appears if one of the sets $I$ or $J$ contains exactly one element, therefore the coefficient of $\dfrac{1}{t_i}$ in the left hand side is equal to $2t_{M\setminus \{i\}}^{|M|-3}$. In the right hand side,
 the factor of $\dfrac{1}{t_i}$ is equal to $2t_M^{|M|-3}$; its value at $t_i=0$ equals $2t_{M\setminus \{i\}}^{|M|-3}$.

\begin{lemma}
The difference
\begin{equation}
\sum_{I \sqcup J = M} t_I^{|I|-2} t_J^{|J|-2} - 2t_M^{|M|-3} \left(\frac{1}{t_1} + \ldots + \frac{1}{t_{|M|}}\right)
\label{diff1}
\end{equation}
is a symmetric polynomial in the variables $t_i$, $i\in M$, and it depends only on the sum of these variables.
\label{l42}
\end{lemma}

The difference is a symmetric function in variables $t_i$, since both summands are symmetric functions in the variables $t_i$. 
Moreover, we know already that the difference in the statement of the Lemma  is a polynomial. 
Therefore, the only thing to prove is that the corresponding symmetric polynomial depends only on $t_M$, the sum of variables $t_i$. It suffices to prove that if we apply the difference operator $\Delta_i$,
$$\Delta_i f(t_1,\ldots,t_i,\ldots,t_m)=f(t_1,\ldots,t_i+h,\ldots,t_m)-f(t_1,\ldots,t_i,\ldots,t_m),$$
then we obtain the same result for each of the variables $t_i,\; i\in M$. Equivalently, we must prove 
that applying the difference $\Delta_i-\Delta_j$ to Eq.~(\ref{diff1}) gives~$0$.
The result of applying $\Delta_i-\Delta_j$ to the sum in the first summand of (\ref{diff1}) is equal to $0$ for all summands corresponding to those subsets $I\subset M$ that contain either both the indexes~$i$ and~$j$, or none of them.
Therefore, we are left with the sum over those partitions of the set $M$ into two nonintersecting subsets
 in which the two indexes belong to different subsets:
\begin{eqnarray*}
(\Delta_i-\Delta_j)\sum_{I\sqcup J= M}t_I^{|I|-2}t_J^{|J|-2}&=&
2(\Delta_i-\Delta_j)\sum_{{I\sqcup J= M\atop I\ni i,J\ni j}}t_I^{|I|-2}t_J^{|J|-2}.
\end{eqnarray*}

Let us look at the right hand side of Eq.~(\ref{diff1}). Make the substitution~$t_i=x$, $t_j=y$. Then the polynomial $t_I^{|I|-2}$, $I\ni i,I\not\ni j$ is taken to the polynomial $(x+t_{I\setminus\{i\}})^{|I\setminus\{i\}|-1}$. Introduce the polynomial $P_{M}(x) = x(x+t_{M})^{|M|-1}$. Now

\begin{eqnarray*}
\sum_{{I\sqcup J= M\atop I\ni i,J\ni j}}t_I^{|I|-2}t_J^{|J|-2}&=&
\sum_{{I\sqcup J= M_{ij}}}\frac{P_{I}(x)}x\frac{P_{J}(y)}y
\end{eqnarray*}
where $M_{ij}=M\setminus\{i,j\}$.
Notice that the subsets $I$ and $J$ are allowed to be empty from now on. 
Notice also that the polynomials $P_M(x)$ are Abel set polynomials of binomial type~\cite{Wis}, 
and hence satisfy the binomial identity:

\begin{equation*}
\sum_{{I\sqcup J= M_{ij}}}\frac{P_{I}(x)}x\frac{P_{J}(y)}y = \frac{P_{M_{ij}}(x+y)}{x y}.
\end{equation*}

Now computation of the result of applying the difference operator $\Delta_i-\Delta_j=\Delta_x-\Delta_y$ is straightforward:

\begin{multline*}
(\Delta_x-\Delta_y)\frac{P_{M_{ij}}(x+y)}{x y} = (x+y+h)(x+y+h+t_{M_{ij}})^{|M_{ij}|-1} \left(\frac{1}{(x+h)y} - \frac{1}{(y+h)x}\right).
\end{multline*}

On the other hand, applying the difference operator $\Delta_i-\Delta_j$ to the second summand in identity (\ref{diff1}) we obtain
\begin{multline*}
(\Delta_i-\Delta_j)(t_M)^{|M|-3}\left(\frac{1}{t_1} + \ldots + \frac{1}{t_{|M|}}\right) =
(t_M+h)^{|M|-3}(\Delta_i-\Delta_j)\left(\frac{1}{t_1} + \ldots + \frac{1}{t_{|M|}}\right)=\\
=(t_M+h)^{|M|-3}\left(\frac{1}{t_i+h}-\frac{1}{t_i} - \frac{1}{t_j+h} +\frac{1}{t_j}\right)  = (t_M+h)^{|M|-3} (t_i+t_j+h)\left(\frac{1}{(t_i+h)t_j}-\frac{1}{(t_j+h)t_i}\right).
\end{multline*}

Now we can complete the proof of Lemma \ref{l42} by noticing that $x=t_i,\; y=t_j$ and $(t_M+h)^{|M|-3} = (x+y+h+t_{M_{ij}})^{|M_{ij}|-1}$.

Therefore, the function
$$
F(t_M) = \sum_{I\sqcup J= M\atop I\neq\emptyset, J\neq\emptyset}t_I^{|I|-2}t_J^{|J|-2} -
t_M^{|M|-4}\left( 2t_M\left(\frac{1}{t_1}+\ldots +\frac{1}{t_{|M|}} \right) -(|M|-2)(|M|-3)\right)
$$
is a symmetric polynomial in the variables $t_i$, $i\in M$, and it depends only on the sum $t_M$ of these variables. It is left to prove that $F(t_M)$ is identically~$0$. It suffices to prove this fact under substitution $t_i=t$ for all $i\in M$,
\begin{multline*}
F(t_M)|_{t_i=t,i=1,\dots,|M|} = \sum\limits_{I\sqcup J= M\atop I\neq\emptyset, J\neq\emptyset} |I|^{|I|-2}t^{|I|-2}|J|^{|J|-2}t^{|J|-2} - (|M|t)^{|M|-4}(|M|^2+5|M|-6) = \\
= t^{|M|-4}\left(\sum\limits_{I\sqcup J= M\atop I\neq\emptyset, J\neq\emptyset} |I|^{|I|-2}|J|^{|J|-2} -  |M|^{|M|-4}(|M|^2+5|M|-6)  \right).
\end{multline*}

Equivalently, we have to prove that

\begin{equation*}
\sum\limits_{i+j=m\atop i,j\neq 0} \binom{m}{i} i^{i-2}j^{j-2} =  m^{m-4}(m^2+5m-6).
\end{equation*}

The last equality follows from the equality of the coefficients at the monomial $x^2y^2$ in the binomial identity for the classical Abel polynomials:

\begin{equation*}
(x+y)(x+y+n)^{n-1} = \sum\limits_{i+j=n} \binom{n}{i}x(x+i)^{i-1}y(y+j)^{j-1}.
\end{equation*}

\renewcommand{\refname}{References}

\end{document}